\documentclass{article}\begin{document}
\centerline{\bf\Large Evaluations of Series of Hyperbolic Functions\normalsize}
\vskip .4in
\centerline{N.D.Bagis}
\centerline{Stenimahou 5 Edessa, Pellas 58200}
\centerline{Edessa, Greece}
\centerline{nikosbagis@hotmail.gr}
\vskip .2in
\[
\]

\textbf{Keywords}: Hyperbolic functions; Series; Evaluations; Elliptic Functions; Special Functions; 
\[
\]
\centerline{\bf Abstract}

\begin{quote}
In this article we give evaluations of certain series of hyperbolic functions, using Jacobi elliptic functions theory. We also define some new functions that enable us to give characterization of not solvable class of series.

\end{quote}

\section{Introduction}

The study of infinite sums and products of certain hyperbolic functions, has been a point of interest in mathematics for more than last 200 years. Works of many great mathematicians such Euler and Gauss firstly, and Jacobi, Eisenstein, Weierstrass, Abel, later, give partial answer to these strange, at first sight, series, by developing the great theory of elliptic functions and modular forms. Later the followers of Jacobi was large in number. A very few of them was, Weber, Ramanujan, Watson, Hardy, Hecke, Poincare, etc.
Nevertheless some of these sums, for example one of them is 
\begin{equation}
\sum^{\infty}_{n=1}\frac{1}{e^{nx}-1}\textrm{, }x>0,
\end{equation}
where not be able to evaluated with the existing theory (we call these class of sums as ''Gost Sums'').\\
Our concern here is to define new functions similar to the classical elliptic functions, having not necessary double periods, and overcome this problem of the evaluation of Gost sums. We also study these new functions and try to find their properties.\\        

Some remarkable evaluations, of closely related sums to (1), are
\begin{equation}
\sum^{\infty}_{n=1}\frac{n^{4\nu+1}}{e^{2\pi n}-1}=\frac{B_{4\nu+2}}{8\nu+4},
\end{equation}
where $\nu$ is positive integer and $B_j$ denotes the $j$th Bernoulli number.\\
Also in [11] it have been proved that
\begin{equation}
\sum_{n\geq0,n-odd}\frac{n^{4\nu+1}}{e^{n\pi}+1}=-\frac{1}{4}Q_{4\nu+1}-\frac{2^{4\nu-1}}{2\nu+1}B_{4\nu+2},
\end{equation}
where $\nu\in\{0,1,2,\ldots\}$ and $Q_{\nu}$ is defined as
\begin{equation}
Q_{\nu}=2\left(\frac{d^{\nu}}{dx^{\nu}}\frac{1}{e^x+1}\right)_{x=0}.
\end{equation}

Continuing for $\nu$ non zero integer the sums
\begin{equation}
\sum^{\infty}_{n=1}\frac{n^{-2\nu-1}}{e^{2\pi n}-1} 
\end{equation}
can evaluated from (see [2]):\\
\\
\textbf{Theorem 1.}\\
Let $a,b>0$ with $ab=\pi^2$, and let $\nu$ be any non zero integer. Then    
$$
a^{-\nu}\left\{\frac{1}{2}\zeta(2\nu+1)+\sum^{\infty}_{n=1}\frac{n^{-2\nu-1}}{e^{2an}-1}\right\}-
(-b)^{-\nu}\left\{\frac{1}{2}\zeta(2\nu+1)+\sum^{\infty}_{n=1}\frac{n^{-2\nu-1}}{e^{2bn}-1}\right\}=
$$
\begin{equation}
=-2^{2\nu}\sum^{\nu+1}_{n=0}(-1)^n\frac{B_{2n}}{(2n)!}\frac{B_{2\nu+2-2n}}{(2\nu+2-2n)!}a^{\nu+1-n}b^n,
\end{equation}
where $\zeta(s)$ is the Riemann zeta function.\\

An example of (6) for $\nu=-1$ is
\begin{equation}
\sum^{\infty}_{n=1}\frac{n}{e^{2n\pi}-1}=\frac{1}{24}-\frac{1}{8\pi}.
\end{equation}
Eisenstein and later Ramanujan consider infinite sums of the form 
\begin{equation}
\sum^{\infty}_{n=1}\frac{n^{2\nu-1}q^n}{1-q^n}\textrm{, }q=e^{-\pi\sqrt{r}}\textrm{, }r>0\textrm{ with  $\nu=1,2,\ldots$}
\end{equation}
and give relations of how one can reduce them and evaluate them, using only the second and third sum of them, (see [2] chapters 14-15). Formulas for evaluating its derivatives was given also by Ramanujan (see [2] chapter 15 Entry 13 and [3]).\\  

The function
\begin{equation}
A=\sum^{\infty}_{n=1}\frac{q^n}{n(1-q^n)}=\sum^{\infty}_{n=1}\frac{1}{n\left(e^{\pi n \sqrt{r}}-1\right)}\textrm{, }q=e^{-\pi\sqrt{r}}\textrm{, }r>0
\end{equation} 
is simply the logarithm of the Ramanujan-Dedekind eta function $f(-q)$, where (see [7] chapters 21,22 and [3])
$$
f(-q)=\prod^{\infty}_{n=1}\left(1-q^n\right)=\exp(-A)=
$$
\begin{equation}
=2^{1/3}\pi^{-1/2}q^{-1/24}(k_r)^{1/12}(k'_r)^{1/3}K(k_r)^{1/2}.
\end{equation}
Hence
\begin{equation}
\frac{dA}{dq}=\frac{e^{2x}}{4}\sum^{\infty}_{n=1}\frac{1}{\sinh^2(nx)}
\end{equation}
and
\begin{equation}
\sum^{\infty}_{n=1}\frac{1}{\sinh^2(nx)}=-4q\frac{d}{dq}\left(\log f(-q)\right)\textrm{, where }q=e^{-2x}\textrm{, }x>0.
\end{equation}
Also (see [11]):
\begin{equation}
\sum^{\infty}_{n=1}\frac{(-1)^nn}{e^{2\pi n/\sqrt{r}}-1}=\frac{1}{8}-\frac{\sqrt{r}}{4\pi}+\frac{rK}{2\pi^2}(E-K).
\end{equation}
Here $K$ and $E$ are the complete elliptic integrals of the first and second kind at modulus $k$, given by $K'/K=\sqrt{r}$. Note also that when $r$ is positive rational number, the $k$ is one of so-called singular moduli and is algebraic number. The elliptic integrals $K,E$ then, can expressed in terms of rational values of the Gamma function, $\pi$ and algebraic numbers (see [9],[10]).

Also in view of [10] we can use elliptic alpha function to get a more relaxed version of (13).\\
It holds
\begin{equation}
\alpha(r)=\frac{\pi}{4K^2}-\sqrt{r}\left(\frac{E}{K}-1\right).
\end{equation}
Solving with respect to $E$ the above formula we get
\begin{equation}
E=\frac{\pi}{4\sqrt{r}K}+K\left(1-\frac{\alpha(r)}{\sqrt{r}}\right).
\end{equation}
Hence (13) becomes
\begin{equation}
2\sum^{\infty}_{n=1}\frac{n}{e^{4\pi n/\sqrt{r}}-1}-\sum_{n\geq 0,n-odd}\frac{n}{e^{2\pi n/\sqrt{r}}-1}=\frac{1}{8}-\frac{\sqrt{r}}{8\pi}-\alpha(r)\frac{\sqrt{r}K^2}{2\pi^2}.
\end{equation}   
Setting $r=4$ in (14) and using (7) along with the fact that $\alpha(4)=6-4\sqrt{2}$, we have
\begin{equation}
\sum_{n\geq0,n-odd}\frac{n}{e^{\pi n}-1}=-\frac{1}{24}+\frac{16\pi}{\Gamma\left(-\frac{1}{4}\right)^4}
\end{equation}
A calculation using\\ 
\\
\textbf{Theorem 2.}(see [14])\\
Let $r>0$ and $q=e^{-\pi\sqrt{r}}$, then
\begin{equation}
1-24\sum^{\infty}_{n=1}\frac{n}{e^{\pi n\sqrt{r}}-1}
=\frac{6}{\pi\sqrt{r}}+\left(1+k^2_r-\frac{6 \alpha(r)}{\sqrt{r}}\right)\frac{4K^2}{\pi^2}.
\end{equation}

can show us that if
\begin{equation}
\mu(r)=-1-9r+12\sqrt{r}\alpha(1/r)+6(r-1)k_r+k_r^2(3r-1),
\end{equation} 
then
\begin{equation}
\sum^{\infty}_{n=1}\frac{n}{e^{4\pi n/\sqrt{r}}-1}=\frac{1}{24}-\frac{\sqrt{r}}{16\pi}+\frac{K'^2}{48\pi^2}\mu(r)\textrm{, }r>0
\end{equation} 
and hence\\
\\
\textbf{Theorem 3.}
\begin{equation}
\sum_{n\geq 0,n-odd}\frac{n}{e^{2\pi n/\sqrt{r}}-1}=-\frac{1}{24}+\frac{\sqrt{r}\alpha(r)K^2}{2\pi^2}+\frac{\mu(r) K'^2}{24\pi^2}
\end{equation} 

Nevertheless Ramanujan gave the evaluation (see [3] last chapter):\\
\\
\textbf{Theorem 4.}(Ramanujan)
\begin{equation}
1-24\sum^{\infty}_{n\geq0,n-odd}\frac{n}{e^{ny}+1}=z^2(1-2x),
\end{equation}
where $y=\pi\sqrt{r}=\pi K'/K$, $z=2K/\pi$ and $x=k^2=k_r^2$ (Ramanujan's notation).\\

Jacobi has given\\
\\
\textbf{Theorem 5.}(Jacobi)\\
For $r>0$ holds
\begin{equation}
\sum_{n\geq0,n-odd}\frac{1}{\cosh\left(n\pi\sqrt{r}/2\right)}=\frac{Kk_r}{\pi}.
\end{equation}

Also from relation (see [13])
\begin{equation}
\frac{\textrm{sn}(q,u)}{\textrm{cn}(q,u)\textrm{dn}(q,u)}=\frac{\pi}{2(k'_r)^2K}\tan\left(\frac{\pi u}{2K}\right)+\frac{2\pi}{(k'_r)^2K}\sum^{\infty}_{n=1}\frac{(-1)^nq^n}{1+q^n}\sin\left(\frac{n\pi u}{K}\right),
\end{equation}
we get\\
\\
\textbf{Theorem 6.}\\
If $r>0$, then
\begin{equation}
\sum^{\infty}_{n=0}\frac{(-1)^n}{e^{(2n+1)\pi\sqrt{r}}+1}=\frac{1}{4}-\frac{Kk'_r}{2\pi}
\end{equation}

\section{Some General Properties and Related Series}

If $x>0$, then under some weak conditions on the sequence $X(n)$ we have the following formula for arithmetical functions $X(n)$:
\begin{equation}
\sum^{\infty}_{n=1}\frac{X(n)}{e^{nx}-1}=\sum^{\infty}_{n=1}e^{-nx}\sum_{d|n}X(d).
\end{equation}  
Hence the Ghost Sum is the generating function of the divisor $d(n)=\sum_{d|n}1$ function i.e.\\
\\
\textbf{Proposition 1.}\\
\begin{equation}
g(x):=\sum^{\infty}_{n=1}\frac{1}{e^{nx}-1}=\sum^{\infty}_{n=1}d(n)q^n\textrm{, where }q=e^{-x}\textrm{, }x>0. 
\end{equation}

This property makes it quite special.\\
Another interesting thing is that ''almost all'' appear in the theory of Jacobian elliptic functions (see [8],[7]). Here we use the Ramanujan's notation in [3] pg.176:\\
\\
\textbf{Proposition 2.}(Ramanujan),([3] pg.174)
\begin{equation}
\sec(\theta)+4\sum^{\infty}_{n=0}\frac{(-1)^n\cos((2n+1)\theta)}{e^{(2n+1)y}-1}=z\sec(\phi)\sqrt{1-x\sin^2(\phi)}.
\end{equation}

Relation (28) is quite close to what we search. 
Setting $\theta=0$ we get (in our notation):\\
\\
\textbf{Corollary 1.}\\
If $q=e^{-\pi\sqrt{r}}$, $r>0$
\begin{equation}
1+4\sum^{\infty}_{n=0}\frac{(-1)^n}{e^{(2n+1)\pi\sqrt{r}}-1}=\frac{2K}{\pi}.
\end{equation}
Hence also
\begin{equation}
1+4\sum_{n\geq0,n\equiv1(4)}\frac{1}{e^{n\pi\sqrt{r}}-1}-4\sum_{n\geq0,n\equiv3(4)}\frac{1}{e^{n\pi\sqrt{r}}-1}=\frac{2K}{\pi}.
\end{equation}
\\
\textbf{Theorem 7.}
\begin{equation}
\sum^{\infty}_{n=1}\frac{n}{\sinh^2(nx)}=-2\frac{d}{dx}\sum^{\infty}_{n=1}\frac{1}{e^{2nx}-1}\textrm{, }x>0.
\end{equation}
\\

The Fourier series of the Jacobi elliptic functions $\textrm{sn}$, $\textrm{cd}$ and $\textrm{cn}$, $\textrm{cd}$, $\textrm{sd}$ are (the functions $\textrm{cn}_1$, $\textrm{cc}$, $\textrm{dd}$, $\textrm{cd}_1$ are not corresponding to elliptic functions and the notation is not the usual. For example $cd=\frac{cn}{dn}$, but $\textrm{dd}$ is not $\frac{dn}{dn}$ so one must be careful):
$$
\textrm{sn}=\textrm{sn}(q,u)=\frac{2\pi}{Kk_r}\sum^{\infty}_{n=0}\frac{q^{n+1/2}\sin\left((2n+1)\frac{\pi}{2K}u\right)}{1-q^{2n+1}}
$$
$$
\textrm{cn}=\textrm{cn}(q,u)=\frac{2\pi}{Kk_r}\sum^{\infty}_{n=0}\frac{q^{n+1/2}\cos\left((2n+1)\frac{\pi}{2K}u\right)}{1+q^{2n+1}}
$$
$$
\textrm{cn}_1=\textrm{cn}_1(q,u)=\frac{2\pi}{Kk_r}\sum^{\infty}_{n=0}\frac{q^{n+1/2}\cos\left((2n+3)\frac{\pi}{2K}u\right)}{1+q^{2n+1}}
$$
$$
\textrm{sd}=\textrm{sd}(q,u)=\frac{2\pi}{Kk_rk'_r}\sum^{\infty}_{n=0}\frac{(-1)^nq^{n+1/2}\sin\left((2n+1)\frac{\pi}{2K}u\right)}{1+q^{2n+1}}
$$
$$
\textrm{cc}=\textrm{cc}(q,u)=\frac{2\pi}{Kk_r}\sum^{\infty}_{n=0}\frac{q^{n+1/2}\cos\left((2n+1)\frac{\pi}{2K}u\right)}{1+q^{2n-1}}=
$$
$$
=\frac{2\pi}{Kk_r}\frac{q^{1/2}\cos(z)}{1+q^{-1}}+q\cdot\textrm{cn}_1(u).
$$
$$
\textrm{cd}=\textrm{cd}(q,u)=\frac{2\pi}{Kk_r}\sum^{\infty}_{n=0}\frac{(-1)^nq^{n+1/2}\cos\left((2n+1)\frac{\pi}{2K}u\right)}{1-q^{2n+1}}
$$
$$
\textrm{dd}=\textrm{dd}(q,u)=\frac{2\pi}{Kk_r}\sum^{\infty}_{n=0}\frac{(-1)^nq^{n+1/2}\cos\left((2n+1)\frac{\pi}{2K}u\right)}{1-q^{2n-1}}
$$
and
\begin{equation}
\textrm{cd}_1=\textrm{cd}_1(q,u)=\frac{2\pi}{Kk_r}\sum^{\infty}_{n=0}\frac{(-1)^nq^{n+1/2}\cos\left((2n+3)\frac{\pi}{2K}u\right)}{1-q^{2n+1}}
\end{equation}
But
$$
\textrm{sn}(-q,u)=\frac{2\pi}{K^{*}k^{*}_r}\sum^{\infty}_{n=0}\frac{i(-1)^nq^{n+1/2}\sin\left((2n+1)\frac{\pi}{2K^{*}}u\right)}{1+q^{2n+1}}=
$$
$$
=i\frac{Kk_rk'_r}{K^{*}k^{*}_r}\frac{2\pi}{Kk_rk'_r}\sum^{\infty}_{n=0}\frac{(-1)^nq^{n+1/2}\sin\left((2n+1)\frac{\pi}{2K}u\frac{K}{K^{*}}\right)}{1+q^{2n+1}}=
$$
$$
=i\frac{Kk_rk'_r}{K^{*}k^{*}_r}\cdot\textrm{sd}\left(q,u\frac{K}{K^{*}}\right)
$$
where we have set $K^{*}=K(k_{r_1})$ and $k^{*}=k_{r_1}$, with $r_1$ such that $e^{-\pi\sqrt{r_1}}=-e^{-\pi\sqrt{r}}$.\\ 
One can see that $k^2_rk^{{*}{2}}_r=k^2_r+k^{{*}{2}}_r$, or the equivalent $k'_r(k^{*}_r)'=1$. Hence from the modular identity
\begin{equation}
\frac{1}{\sqrt{1-x}}K\left(\frac{x}{x-1}\right)=K(x),
\end{equation}
we get
\begin{equation}
\frac{K(k_{r_1})}{K(k_{r})}=\frac{K^{*}}{K}=k'_r.
\end{equation}
\\
\textbf{Proposition 3.}\\
If $q=e^{-\pi\sqrt{r}}$, $r>0$, then
\begin{equation}
\textrm{sn}(-q,u)=k'_r\cdot\textrm{sd}\left(q,\frac{u}{k'_r}\right).
\end{equation}

Consider now Theorem 5 in the form
\begin{equation}
2\sum_{n\geq0,n-odd}\frac{q^{n/2}}{1+q^n}=\frac{Kk_r}{\pi}
\end{equation}
and set $q\rightarrow -q$, then from the above relations (34) and $k^{*}_r=ik_r/k'_r$, we get\\
\\
\textbf{Theorem 8.}\\
If $q=e^{-\pi\sqrt{r}}$, $r>0$, then
\begin{equation}
2\sum^{\infty}_{n\geq0,n-odd}\frac{(-1)^{\frac{n-1}{2}}q^{n/2}}{1-q^n}=\sum^{\infty}_{n=0}\frac{(-1)^{n}}{\sinh\left((n+1/2)\pi\sqrt{r}\right)}=\frac{Kk_r}{\pi}.
\end{equation}

If 
\begin{equation}
\textrm{ss}:=\textrm{ss}(q,u):=\frac{2\pi}{Kk_r}\sum^{\infty}_{n=0}\frac{q^{n+1/2}\sin\left((2n+1)\frac{\pi}{2K}u\right)}{1+q^{2n+1}},
\end{equation}
then from the above formula of $\textrm{cn}$ and the elementary trigonometric formula\\ $\cos(a+b)=\cos(a)\cos(b)-\sin(a)\sin(b)$ we get  
$$
\textrm{cc}(q,u)=\frac{2\pi}{Kk_r}\frac{q^{1/2}\cos(z)}{1+q^{-1}}+q\cdot\textrm{cn}_1(q,u)=\frac{2\pi}{Kk_r}\frac{q^{1/2}\cos(z)}{1+q^{-1}}+
$$
$$
+\frac{2\pi}{Kk_r}\sum^{\infty}_{n=0}\frac{q^{n+1/2+1}}{1+q^{2n+1}}\left(\cos((2n+1)z)\cos(2z)-\sin((2n+1)z)\sin(2z)\right)=
$$
$$
=\frac{2\pi}{Kk_r}\frac{q^{1/2}\cos(z)}{1+q^{-1}}+q\cos(2z) \textrm{cn}(q,u)-q\sin(2z)\textrm{ss}(q,u).
$$
Hence we get the next\\
\\
\textbf{Proposition 4.}\\
If $q=e^{-\pi\sqrt{r}}$, $r>0$
$$
\textrm{ss}(q,u)=\cot(2z)\textrm{cn}(q,u)-q^{-1}\csc(2z)\textrm{cc}(q,u)+\frac{\pi\sqrt{q}\csc(z)}{(1+q)k_rK}=
$$
\begin{equation}
=\textrm{cn}(q,u)\cot(2z)-\textrm{cn}_1(q,u)\csc(2z),
\end{equation}
where $z=\frac{\pi u}{2K}$.\\
\\ 
\textbf{Proposition 5.}\\
When $q=e^{-\pi\sqrt{r}}$, $r>0$, it holds
\begin{equation}
\textrm{cn}(-q,u)=\textrm{cd}\left(q,\frac{u}{k'_r}\right)
\end{equation}
\\
\textbf{Proof.}\\
From the $\textrm{cn}$ formula we have
$$
\textrm{cn}(-q,u)=\frac{2\pi}{K^{*}k^{*}_r}\sum^{\infty}_{n=0}\frac{(-1)^n q^{n+1/2}}{1-q^{2n+1}}\cos\left((2n+1)\frac{\pi u}{2K^{*}}\right)=
$$
$$
=i\frac{Kk_r}{K^{*}k^{*}_r}\frac{2\pi}{Kk_r}\sum^{\infty}_{n=0}\frac{(-1)^n q^{n+1/2}}{1-q^{2n+1}}\cos\left((2n+1)\frac{\pi u}{2K}\frac{K}{K^{*}}\right)=
$$
$$
=i\frac{Kk_r}{K^{*}k^{*}_r}\textrm{cd}\left(q,u\frac{K}{K^{*}}\right)=\textrm{cd}\left(q,\frac{u}{k'_r}\right).
$$
Since 
\begin{equation}
i\frac{Kk_r}{K^{*}k^{*}_r}=1.
\end{equation}
$qed$.

\section{The Function $\textrm{cd}_1$ : Evaluations and Properties}

\textbf{Theorem 9.}\\
If $a=iq^{1/2}e^{i\pi\theta/(2K)}$ and $\theta\in \textbf{C}$, with $|q^{1/2}e^{i\pi\theta/(2K)}|<1$, then
\begin{equation}
\frac{2\pi }{Kk_r}\sum^{\infty}_{n=0}\frac{a^{2n+1}}{1-q^{2n+1}}=-\textrm{cd}\left(q,\theta\right)\cot\left(\frac{\theta\pi}{K}\right)+\textrm{cd}_1\left(q,\theta\right)\csc\left(\frac{\theta\pi}{K}\right)+i\cdot\textrm{cd}\left(q,\theta\right)
\end{equation}
\\
\textbf{Proof.}\\
It is easy to see someone that
$$
\frac{2\pi }{Kk_r}\sum^{\infty}_{n=0}\frac{a^{2n+1}}{1-q^{2n+1}}=
i\cdot\textrm{cn}(-q,\theta k'_r)-\textrm{ss}\left(-q,\theta k'_r\right)=
$$
$$
=i\cdot \textrm{cd}\left(q,\theta\right)-\textrm{ss}\left(-q,\theta k'_r\right)
$$
But from Propositions 4,5 we have
$$
\textrm{ss}(-q,u)=\textrm{cd}\left(q,\frac{u}{k'_r}\right)\cot(2z^{*})-\textrm{cn}_1(-q,u)\csc(2z^{*})\textrm{, }z^{*}=\frac{\pi u}{2K^{*}}.
$$
Also
\begin{equation}
\textrm{cn}_1(-q,u)=\textrm{cd}_1\left(q,\frac{u}{k'_r}\right)
\end{equation}
Hence
\begin{equation}
\textrm{ss}(-q,u)=\textrm{cd}\left(q,\frac{u}{k'_r}\right)\cot\left(\frac{2z}{k'_r}\right)-\textrm{cd}_1\left(q,\frac{u}{k'_r}\right)\csc\left(\frac{2z}{k'_r}\right).
\end{equation}
From the above relations we get the proof. $qed$\\
\\
\textbf{Theorem 10.}\\
If $0<|\lambda|<1$, $r>0$, then
$$
\frac{\pi}{Kk_r}\sum^{\infty}_{n=0}\frac{(-1)^ne^{-\pi\sqrt{r}(n+1/2)\lambda}}{\sinh\left((n+1/2)\pi\sqrt{r}\right)}=\textrm{cd}\left(q,\lambda i K'\right)\coth\left(\lambda\pi\sqrt{r}\right)-
$$
\begin{equation}
-\textrm{cd}_1\left(q,\lambda i K'\right)\textrm{csch}\left(\lambda\pi\sqrt{r}\right)+\textrm{cd}\left(q,\lambda i K'\right).
\end{equation}
\\
\textbf{Proof.}\\
Set $\theta=\lambda i K'$ in relation (42) of Theorem 9. $qed$\\
\\
\textbf{Theorem 11.}\\
If $q=e^{-\pi\sqrt{r}}$, $r>0$ and $\nu\in \textbf{C}^{*}$, such that $2/\nu$ not integer, then
$$
\frac{2\pi}{Kk_r}\sum^{\infty}_{n=0}\frac{q^{(2n+1)(1/2+1/\nu)}}{1-q^{2n+1}}=i\cdot\textrm{sn}\left(q,\frac{2iK'}{\nu}\right)\coth\left(\frac{2\pi\sqrt{r}}{\nu}\right)+
$$
\begin{equation}
+i\cdot\textrm{cd}_1\left(q,-K+\frac{2iK'}{\nu}\right)\textrm{csch}\left(\frac{2\pi\sqrt{r}}{\nu}\right)+i\cdot\textrm{sn}\left(q,\frac{2iK'}{\nu}\right).
\end{equation}
In case that $r\in\textbf{Q}^{*}_{+}$ and $\nu\in\textbf{Q}^{*}_{+}-\textbf{Z}$, then $\textrm{sn}\left(q,\frac{2iK'}{\nu}\right)$ is algebraic number.\\ 
\\
\textbf{Proof.}\\
If we replace 
$$
\theta=\theta_1:=-i\log\left(-iq^{1/\nu}\right)\frac{2K}{\pi}=\left(-\frac{\pi}{2}+i\frac{\pi\sqrt{r}}{\nu}\right)\frac{2K}{\pi}=
$$
\begin{equation}
=\left(-1+i\frac{2\sqrt{r}}{\nu}\right)K=-K+\frac{2iK'}{\nu},
\end{equation}
then we will have $a=q^{1/2+1/\nu}$. From the relations (see [8]):
\begin{equation}
\textrm{cn}\left(q,u+K\right)=-k'_r\textrm{sn}(q,u)/\textrm{dn}(q,u)
\end{equation}
and
\begin{equation}
\textrm{dn}\left(q,u+K\right)=k'_r/\textrm{dn}(q,u)
\end{equation}
we have
\begin{equation}
\textrm{cd}\left(q,u+K\right)=-\textrm{sn}(q,u)
\end{equation}
and easy
$$
\textrm{cd}\left(q,u-K\right)=-\textrm{sn}\left(q,u-2K\right)=\textrm{sn}(q,u),
$$
since $\textrm{sn}\left(q,u+2K\right)=-\textrm{sn}(q,u)$.\\Also easily we get
$$
\csc\left(\frac{\theta_1\pi}{K}\right)=i\cdot\textrm{csch}\left(\frac{2\pi\sqrt{r}}{\nu}\right)$$ 
and 
$$
\cot\left(\frac{\theta_1\pi}{K}\right)=-i\cdot\textrm{coth}\left(\frac{2\pi\sqrt{r}}{\nu}\right).
$$
From the above and Theorem 9 we get the result. $qed$\\ 
\\
\textbf{Corollary 2.}\\
If $q=e^{-\pi\sqrt{r}}$, $r>0$, then\\
1)
\begin{equation}
\lim_{y\rightarrow K}\frac{\textrm{cd}_1\left(q,y\right)}{y-K}=1+\frac{2\pi^2}{K^2k_r}\sum_{n\geq0\textrm{, }n-odd}\frac{q^{n/2}}{1-q^n},
\end{equation}
2) $\textrm{cd}_1(q,0)=1$, $\textrm{cd}_1(q,K)=0$, $\textrm{cd}_1(q,2K)=-1$ and $\textrm{cd}_1(u+2K)=-\textrm{cd}_1(u)$.\\
\\
\textbf{Proof.}\\
Taking the limit $\nu\rightarrow\infty$ in (46) and using
$$
\lim_{\nu\rightarrow\infty}\textrm{sn}\left(q,\frac{2iK'}{\nu}\right)=0,
$$
$$
\lim_{\nu\rightarrow\infty}\textrm{sn}\left(q,\frac{2iK'}{\nu}\right)\coth\left(\frac{2\pi\sqrt{r}}{\nu}\right)=\frac{iK'}{\pi\sqrt{r}},
$$
we get easily the result. $qed$\\
\\
\textbf{Theorem 12.}\\
If $q=e^{-\pi\sqrt{r}}$, $r>0$, then
\begin{equation}
\textrm{cd}_1(q,iK')=\frac{1}{qk_r}-\frac{\sinh(\pi\sqrt{r})}{k_r}\left(1-\frac{\pi}{2K}\right).
\end{equation}
\\
\textbf{Proof.}\\
Set $\theta=iK'$ in (42). Then using relation $\textrm{cd}\left(q,iK'\right)=\frac{1}{k_r}$ (see [8]), we get 
$$
\textrm{cd}_1\left(q,iK'\right)=\frac{1}{qk_r}-\frac{2\pi\sinh\left(\pi\sqrt{r}\right)}{Kk_r}\sum^{\infty}_{n=0}\frac{(-1)^n}{e^{(2n+1)\pi\sqrt{r}}-1}.
$$
The result follows from Corollary 1. $qed$\\
\\
\textbf{Notes.}\\
i) Numerical values of $\textrm{cd}_1(q,iK')$ can given using Theorems 13 and 14 bellow.\\ 
ii) Formula (32) does not converges for these values.\\
\\
\textbf{Corollary 3.}\\
If $q=e^{-\pi\sqrt{r}}$, $r>0$, then
\begin{equation}
\textrm{cd}_1\left(q,\frac{iK'}{2}\right)=\frac{1}{\sqrt{qk_r}}-\frac{\pi\sinh\left(\frac{\pi\sqrt{r}}{2}\right)}{Kk_r}\sum^{\infty}_{n=0}\frac{(-1)^ne^{-(n+1/2)\pi\sqrt{r}/2}}{\sinh\left((n+1/2)\pi\sqrt{r}\right)}
\end{equation}
\\
\textbf{Corollary 4.}\\
If $q=e^{-\pi\sqrt{r}}$, $r>0$, then
\begin{equation}
\textrm{cd}_1\left(q,\frac{K}{2}\right)=\frac{1}{\sqrt{1+k'_r}}-\frac{\pi\sqrt{8}}{Kk_r}\sum_{n\geq0\textrm{,  }n\equiv 1(4)}\frac{(-1)^{\frac{n-1}{4}}q^{n/2}}{1-q^n}
\end{equation}
and
\begin{equation}
-\sum^{\infty}_{n=1}\chi(n)\frac{q^{n/2}}{1-q^n}=\frac{Kk_r}{\sqrt{2}\pi\sqrt{1+k'_r}},
\end{equation}
where $\chi(n)=\left(\frac{n+2}{8}\right)$ and $\left(\frac{n}{m}\right)$ is the usual Jacobi symbol.\\ 
\\
\textbf{Proof.}\\
Setting $\theta=\frac{K}{2}$ in Theorem 9, we get $a=iq^{1/2}e^{\pi i/4}=q^{1/2}e^{3\pi i/4}$, $\cot\left(\frac{\pi}{2}\right)=0$, $\csc\left(\frac{\pi}{2}\right)=1$. Also it holds (see [8]):
\begin{equation}
\textrm{cn}\left(q,\frac{K}{2}\right)=\frac{\sqrt{k'_r}}{\sqrt{1+k'_r}}\textrm{ and  }\textrm{dn}\left(q,\frac{K}{2}\right)=\sqrt{k'_r}.
\end{equation}
Hence
\begin{equation}
\textrm{cd}\left(q,\frac{K}{2}\right)=\frac{\textrm{cn}\left(q,\frac{K}{2}\right)}{\textrm{dn}\left(q,\frac{K}{2}\right)}=\frac{1}{\sqrt{1+k'_r}}.
\end{equation}
From the above we can evaluate
\begin{equation}
\textrm{cd}_1\left(q,\frac{K}{2}\right)=e^{3\pi i/4}\frac{2\pi}{Kk_r}\sum^{\infty}_{n=0}\frac{q^{n+1/2}e^{3\pi i n/2}}{1-q^{2n+1}}-i\cdot\textrm{cd}\left(q,\frac{K}{2}\right).
\end{equation}
Hence taking the real and imaginary parts of the above equation we deduce 
\begin{equation}
\sum_{n\geq 0\textrm{, }n\equiv 1(4)}\frac{(-1)^{\frac{n-1}{4}}q^{n/2}}{1-q^n}+\sum_{n\geq 0\textrm{, }n\equiv 3(4)}\frac{(-1)^{\frac{n-3}{4}}q^{n/2}}{1-q^n}=\frac{\sqrt{2}}{\sqrt{1+k'_r}}\frac{Kk_r}{2\pi}
\end{equation}
and
\begin{equation}
\frac{\pi\sqrt{2}}{Kk_r}\sum_{n\geq 0\textrm{, }n\equiv 3(4)}\frac{(-1)^{\frac{n-3}{4}}q^{n/2}}{1-q^n}-\frac{\pi\sqrt{2}}{Kk_r}\sum_{n\geq 0\textrm{, }n\equiv 1(4)}\frac{(-1)^{\frac{n-1}{4}}q^{n/2}}{1-q^n}=\textrm{cd}_1\left(q,\frac{K}{2}\right).
\end{equation}
$qed$\\

Continuing we define 
\begin{equation}
(z;q)_{\infty}:=\prod^{\infty}_{n=0}(1-zq^n)\textrm{, }|q|<1\textrm{, }z\in\textbf{C}.
\end{equation}
also
\begin{equation}
P:=\left(\frac{(-a;q)_{\infty}}{(a;q)_{\infty}}\right)^2
\end{equation}
and
\begin{equation}
u_0(a,q):=\frac{P-1}{P+1}
\end{equation}
Then we have (see [12]):
\begin{equation}
\log\left(-1+\frac{2}{1-u_0(a,q)}\right)=\log P
\end{equation}
and
\begin{equation} \log P=\log\left(-1+\frac{2}{1-u_0(a,q)}\right)=4\sum^{\infty}_{n=0}\frac{a^{2n+1}}{(2n+1)(1-q^{2n+1})}.
\end{equation}
If $\theta$ is a complex number, we set $a=iq^{1/2}e^{i\pi t/(2K)}$ and derivate (65) with respect to parameter $t$, then set the value $t=\theta$.
$$
\left[\frac{d}{dt}\log\left(-1+\frac{2}{1-u_0(a,q)}\right)\right]_{t=\theta}=4\left[\frac{d}{dt}\sum^{\infty}_{n=0}\frac{a^{2n+1}}{(2n+1)(1-q^{2n+1})}\right]_{t=\theta}=
$$
\begin{equation}
=\frac{2\pi i}{K}\sum^{\infty}_{n=0}\frac{a^{2n+1}}{1-q^{2n+1}}.
\end{equation}
Hence using Theorem 9 and the next integral relation [8]
\begin{equation}
\int \textrm{cd}(q,t)dt=\log\left(\textrm{nd}(q,t)+k_r\textrm{sd}(q,t)\right),
\end{equation}
we get\\
\\
\textbf{Theorem 13.}\\
If $q=e^{-\pi\sqrt{r}}$, $r>0$ and $a=iq^{1/2}e^{i\pi \theta/(2K)}$, $\theta$ real number, then 
\begin{equation}
Re\left[\log\left(-1+\frac{2}{1-u_0(a,q)}\right)\right]
=\log\left(\textrm{nd}(q,\theta)+k_r\textrm{sd}(q,\theta)\right)
\end{equation}
and
\begin{equation}
\textrm{cd}_1(q,\theta)=\textrm{cd}(q,\theta)\cos\left(\frac{\pi\theta}{K}\right)+2k_r^{-1}\sin\left(\frac{\pi\theta}{K}\right)Im\left[\frac{d}{d\theta}\log\left(\frac{\left(-a;q\right)_{\infty}}{\left(a;q\right)_{\infty}}\right)\right].
\end{equation}
\\
\textbf{Note.}\\
We have the next continued fraction expansion (see [12]):
\begin{equation}
u_0(a,q)=\frac{2a}{1-q+}\frac{a^2(1+q)^2}{1-q^3+}\frac{a^2q(1+q^2)^2}{1-q^5+}\frac{a^2q^2(1+q^3)^2}{1-q^7+}\ldots\textrm{, }|q|<1.
\end{equation}
This continued fraction can be used to get numerical verifications of values such $\textrm{cd}_1(q,i\cdot nK')$, $n\in\textbf{Q}^{*}_{+}$.\\ 
\\
\textbf{Theorem 14.}\\
If $q=e^{-\pi\sqrt{r}}$, $r>0$ and $a=iq^{1/2}e^{i\pi t/(2K)}$, $t$ parameter and $\theta\in\textbf{C}$ such that $|q^{1/2}e^{i\pi\theta/(2K)}|<1$, then
$$
\textrm{cd}_1(q,\theta)=\textrm{cd}(q,\theta)\cos\left(\frac{\pi\theta}{K}\right)-i\cdot\textrm{cd}(q,\theta)\sin\left(\frac{\pi\theta}{K}\right)-i\cdot k_r^{-1}\sin\left(\frac{\pi\theta}{K}\right)\times
$$
\begin{equation}
\times\left[\frac{d}{dt}\log\left(-1+\frac{2}{1+}\frac{-2a}{1-q+}\frac{a^2(1+q)^2}{1-q^3+}\frac{a^2q(1+q^2)^2}{1-q^5+}\frac{a^2q^2(1+q^3)^2}{1-q^7+}\ldots\right)\right]_{t=\theta}.
\end{equation}
and
$$
\textrm{cd}_1(q,\theta)=\textrm{cd}(q,\theta)\cos\left(\frac{\pi\theta}{K}\right)-i\cdot\textrm{cd}(q,\theta)\sin\left(\frac{\pi\theta}{K}\right)
-
$$
\begin{equation}
-2i\cdot k_r^{-1}\sin\left(\frac{\pi\theta}{K}\right)\left[\frac{d}{dt}\log\left(\frac{(-a;q)_{\infty}}{(a;q)_{\infty}}\right)\right]_{t=\theta}.
\end{equation}
\\
\textbf{Proof.}\\
Using (62),(65),(66),(70),(42), we get the two results. $qed$\\
\\
\textbf{Theorem 15.}\\
Let $1/\nu$ be positive integer. Let also $m$ be even integer and $\nu_1=2/\nu$, then if $q=e^{-\pi\sqrt{r}}$, $r>0$ we have
$$
\textrm{cd}_1\left(q,mK+\nu_1iK'\right)=(-1)^{m/2}e^{\nu_1\pi \sqrt{r}}-
$$
\begin{equation}
-(-1)^{m/2}\frac{2\pi }{Kk_r}\sinh\left(\nu_1\pi\sqrt{r}\right)\left(-\sum^{1/\nu-1}_{j=0}\frac{q^{j+1/2}}{1+q^{2j+1}}
+\frac{Kk_r}{2\pi}\right).
\end{equation} 
\\
\textbf{Proof.}\\
The proof follows from the identities (66),(71) along with
\begin{equation}
\sum^{\infty}_{n=0}\frac{(-1)^nq^{(2n+1)(l+1/2)}}{1-q^{2n+1}}=-\sum^{l-1}_{j=0}\frac{q^{j+1/2}}{1+q^{2j+1}}+\frac{Kk_r}{2\pi}\textrm{, }l\in\textbf{N}\textrm{, }q=e^{-\pi\sqrt{r}}\textrm{, }r>0,
\end{equation}
and 
\begin{equation}
\textrm{cd}(q,mK+\nu_1K')=(-1)^{m/2},
\end{equation}
with $m$ even integer and $\nu_1=2/\nu$ positive integer. $qed$

\[
\]

\centerline{\bf References}\vskip .2in

\noindent

[1]: M.Abramowitz and I.A.Stegun, 'Handbook of Mathematical Functions'. Dover Publications, New York., (1972).

[2]: B.C. Berndt, 'Ramanujan`s Notebooks Part II'. Springer Verlang, New York., (1989).

[3]: B.C. Berndt, 'Ramanujan`s Notebooks Part III'. Springer Verlang, New York., (1991).

[4]: I.S. Gradshteyn and I.M. Ryzhik, 'Table of Integrals, Series and Products'. Academic Press., (1980).

[5]: L. Lorentzen and H. Waadeland, Continued Fractions with Applications. Elsevier Science Publishers B.V., North Holland., (1992).  

[6]: H.S. Wall. 'Analytic Theory of Continued Fractions'. Chelsea Publishing Company, Bronx, N.Y., (1948). 

[7]: E.T. Whittaker and G.N. Watson. 'A course on Modern Analysis'. Cambridge U.P., (1927).

[8]: J.V. Armitage, W.F. Eberlein. 'Elliptic Functions'. Cambridge University Press., (2006).

[9]: J.M. Borwein, M.L. Glasser, R.C. McPhedran, J.G. Wan, I.J. Zucker. 'Lattice Sums Then and Now'. Cambridge University Press. New York., (2013).

[10]: J.M. Borwein and P.B. Borwein. 'Pi and the AGM: A Study in Analytic Number Theory and Computational Complexity', Wiley, New York., (1987).

[11]: M.L. Glasser and N.D. Bagis. 'Some Applications of the Poisson Summation Formula'. arXiv:0812.0990, (2008)  

[12]: N.D. Bagis and M.L. Glasser. 'Evaluations of a Continued Fraction of Ramanujan'. Rend. Sem. Mat. Univ. Padova. Vol 133., (2015). 

[13]: S.C. Milne. 'Infinite Families of Exact Sums of Squares Formulas, Jacobi Elliptic Functions, Continued Fractions, and Schur Functions'.\\ arXiv:math/0008068v2 [math.NT] 7 Juan. 2001.

[14]: N.D. Bagis and M.L. Glasser. 'On the Transcendence of Complete Elliptic Integrals of the First Kind and Values of the Gamma Function'. submitted

\end{document}